\documentclass{article}
\usepackage[dvips]{graphicx}
\usepackage{amsmath}
\usepackage{amssymb}
\usepackage{hhline}
\usepackage{amscd}

\renewcommand\a{\alpha}
\renewcommand\b{\beta}

\renewcommand\c{\circ}
\newcommand\C{{\mathbb C}}
\newcommand\CC{{\cal C}}

\renewcommand\d{\delta}

\newcommand\e{\eta}

\newcommand\g{\gamma}
\newcommand\G{\Gamma}

\newcommand\hps{\hskip-16pt . \hskip2pt}

\newcommand\J{{\hat J}}

\renewcommand\l{\lambda}

\newcommand\La{\Lambda}

\newcommand\na{\nabla}
\renewcommand\O{\Omega}
\newcommand\oo{\omega}
\newcommand\ot{\otimes}
\newcommand\ok{{\cal O}}

\newcommand\op[1]{\mathop{\rm #1}\nolimits}
\newcommand\p{\partial}
\newcommand\po{$\!\!\!{\bf .}$ }
\renewcommand\P{\Phi}

\newcommand\R{{\mathbb R}}

\newcommand\te{\theta}
\newcommand\ti{\tilde}

\newcommand\ve{\varepsilon}
\newcommand\vp{\varphi}

\newcommand\vph{\vphantom{$\frac{2^2}{2_2}$}}

\newcommand\we{\wedge}
\newcommand\x{\xi}

\newcommand\z{\sigma}
\newcommand\Z{{\mathbb Z}}

\def\Rom#1{\uppercase\expandafter{\romannumeral#1}}

\newcommand\1{{\bf 1}}

\newcommand\qed{\phantom{\underline{y}}\hfill\hfill$\square$}

\newcommand\bib[1]{\bibitem[#1]{#1}}

\newcommand\const{\mathop{const}\nolimits}

 \newcommand{\ddfrac}[2]{\dfrac{\mathstrut\displaystyle #1}
 {\mathstrut\displaystyle #2}}

\newtheorem{theorem}{Theorem}
\newtheorem{thea}{Theorem}

\newtheorem{theu}{Theorem}

\newtheorem{cor}[theorem]{Corollary}
\newtheorem{dfn}[thea]{Definition}
\newtheorem{lem}[theorem]{Lemma}
\newtheorem{prop}[theorem]{Proposition}
\newtheorem{rk}[theu]{Remark}
\newenvironment{proof}[1][Proof]{\textbf{#1.} }{\qed}

\newenvironment{ex}{\trivlist \item[\hskip \labelsep{\bf Example.}]}{\endtrivlist}

\begin{document}

\title{{\bf Tangent and normal bundles \\
in almost complex geometry}}
\author{Boris~S.\ Kruglikov}
 \date{}

\maketitle

 \begin{abstract}
We define and study pseudoholomorphic vector bundles structures,
particular cases of which are tangent and normal bundle almost
complex structures. As an application we deduce normal forms of
1-jets of almost complex structures along a pseudoholomorphic
submanifold. In dimension four we relate these normal forms to the
problem of pseudoholomorphic foliation of a neighborhood of a
curve and the question of non-deformation and persistence of
pseudoholomorphic tori.\footnote{MSC-2000 numbers: 53C15, 53A55,
32G05; 53C05, 58D27, 58A20.
 }

 \end{abstract}


\section*{Introduction}

 \hspace{13.5pt}
In this paper we study the differential geometry of tangent and
normal bundles in the almost complex category. Let $J:TM\to TM$ be
an almost complex structure, $J^2=-\1$. A submanifold $L\subset M$
is called {\it pseudoholomorphic\/} (PH-submanifold) if $TL\subset
TM$ is $J$-invariant.

We introduce two different canonical almost complex structures
$\hat J$ and $\check J$ on each of the total spaces $TL$ and
$N_LM$ of tangent and normal bundles such that the projection to
$L$ and the zero section embeddings of $L$ are pseudoholomorphic.
We find an explicit relation between these two almost complex
structures.

Moreover, we define and investigate the theory of abstract
pseudoholomorphic (almost holomorphic) vector bundles, partial
cases of which are tangent and normal bundles. We describe their
normal forms, which produce normal forms of 1-jets of almost
complex structures along a PH-submanifolds.

Generically the only PH-submanifolds are PH-curves (\cite{K2}).
Local existence of PH-curves was established by Nijenhuis and
Woolf (\cite{NW}). The global existence result is due to Gromov,
whose paper \cite{Gro} made compact PH-curves an indispensable
tool of symplectic geometry.

For a PH-curve $L$ the structure $\hat J$ on $N_LM$ is
holomorphic, while the structure $\check J$ is not, and they both
play an important role in the deformation and regularity questions
for PH-curves. In particular, we relate Gromov's operator $D_u$ to
our normal bundle structures. Consequently, the structure $\check
J$ appears to be basic for local Gromov-Witten theory.

In~\cite{Mo} Moser constructed a KAM-type theory for a
PH-foliation of an almost complex torus $T^{2n}$ by non-compact
curves, namely entire PH-lines $\C\to T^{2n}$ with generic slope.
He proved that under a small almost complex perturbation of the
standard complex structure $J_0$ many leaves persist. If the
perturbation is big, but tame-restricted, then only some of the
leaves persist. This was proven by Bangert in~\cite{B}. Another
proof is given in~\cite{KO}.

In \cite{A2}(1993-25) Arnold asks about almost complex version, in
the spirit of Moser's result, for his Floquet-type theory of
elliptic curves neighborhoods (\cite{A1}). It will be shown the
direct extension fails (there are moduli in normal forms), though
we conjecture the right generalization, treatable by Moser's
method, is a possibility of foliation of a PH-torus neighborhood
by PH-cylinders.

We consider specially the case $\dim M=4$ and find the condition
for a PH-curve neighborhood to admit a PH-foliation of a special
kind. We also study problems of persistence and isolation of
PH-tori, as posed by Moser. In particular, we obtain a geometric
interpretation for his non-deformable example from \cite{Mo}.
There Moser announced "a study of the normal bundle", which has
not been performed. The present paper fills the gap.

In appendix A we give a new proof of a theorem by Lichnerowicz,
essentially used in the main constructions, and consider some
applications of the minimal connections. In appendix B we discuss
what happens with the normal and tangent bundles for other
geometric structures, which demonstrates, in particular, that
relations (\ref{apa})-(\ref{ape}) from \S\ref{sec2} is a PH-analog
of the Ricci equation.

\section{\hps Almost complex tangent bundle}
\label{sec1}

 \hspace{13.5pt}
The Nijenhuis tensor of an almost complex structure $J\in
C^\infty(T^*M\ot TM)$ is given by the formula
 $$
N_J(X,Y)=[JX,JY]-J[JX,Y]-J[X,JY]-[X,Y],\quad X,Y\in TM.
 $$
We write $N_J\in C^\infty(\La^2T^*M\ot_{\bar\C}TM)$ meaning it is
skew-symmetric in $X,Y$ and $J$-antilinear. By the
Newlander-Nirenberg theorem \cite{NW} integrability of $J$ can be
expressed as $N_J=0$.

An almost complex connection is a linear connection $\na$ that
preserves the almost complex structure: $\na J=0$. It is called
minimal if its torsion $T_\na=\frac14N_J$. Such connections always
exist due to \cite{L}, see appendix~\ref{a_A}.

Let $\pi:TM\to M$ denote the projection and $\rho:M\to TM$ the
zero section.

 \begin{theorem}\po\label{thm1}
There exists a canonical almost complex structure $\hat J$ on the
total space of the tangent bundle $TM$\,%
to an almost complex manifold $(M,J)$ such that:
 \begin{enumerate}
  \item
The maps $\pi:TM\to M$ and $\rho:M\to TM$ are pseudoholomorphic.
  \item
$(TM,\hat J)$ is integrable iff $(M,J)$ is integrable.
   \end{enumerate}
 \end{theorem}

 \begin{proof}
Consider a minimal connection $\na$. It produces the splitting
$T_a(TM)=H_a\oplus V_a$ into horizontal and vertical components,
$a\in TM$. We have natural isomorphisms $\pi_*:H_a\simeq T_xM$ and
$V_a\simeq T_xM$, $x=\pi(a)$. Define the structure $\hat J$ on
$T_a(TM)$ as $J\oplus J$ with respect to the above splitting and
isomorphisms.

If we change the minimal connection $\tilde\nabla=\nabla+A$
(Theorem~\ref{th_A} of appendix \ref{a_A}), then
$A\in C^\infty(S^2T^*M\ot_\C TM)$. The new horizontal space is
given by ${\ti H}_a=\op{graph}\{A(a,\cdot): H_a\to V_a\}$. Since
$A(a,\cdot)$ is a complex linear map, the almost complex structure
$\hat J$ on $TM$ is defined canonically.

The properties of $\hat J$ follow directly from the construction.
 \end{proof}

 \begin{rk}\po
Whenever integrable, $\hat J$ defines the standard holomorphic
structure.
 \end{rk}

Construction of the structure $\hat J$ can be generalized to the
cotangent and other tensor bundles. The adjoint $J^*$ to the
operator $J$ is a fiberwise complex structure on $T^*M$. The two
structures induce a canonical fiberwise complex structure on the
complex-linear tensor bundles $T^{(r,s)}_{\C}M$ of contravariant
degree $r$ and covariant degree $s$ tensors and also on the
subbundles $S^k_\C TM$, $\La^k_\C TM$. As usual, the tensor
product over $\C$ is formed by the equivalence relation $X\ot
JY\sim JX\ot Y$ (so that $T^{(r,s)}_{\C}M\ne T^{(r,s)}M\ot\C$
etc).

 \begin{theorem}\po\label{thm1,5}
Let $E_M$ be one of the bundles $T^{(r,s)}_{\C}M$, $S^k_\C M$,
$\La^k_\C M$ or their duals and tensor products over $\C$. There
exists a canonical almost complex structure $\hat J$ on the total
space $E_M$ such that:
 \begin{enumerate}
  \item
The maps $\pi:E_M\to M$ and $\rho:M\to E_M$ are pseudoholomorphic.
  \item
$(E_M,\hat J)$ is integrable iff $(M,J)$ is integrable.
   \end{enumerate}
 \end{theorem}

 \begin{proof}
The claim is obtained similarly to Theorem \ref{thm1} by checking
that the admissible gauge transformations
$A^E\in\O^1(M,\op{end}_\C E_M)$ are complex-linear in all
arguments. This follows from the explicit formulae:
$A^{(1,0)}(X)=A(X)$, $A^{(0,1)}(X)=-A(X)^*$,
$A^{(2,0)}(X)=A(X)\ot\1+\1\ot A(X)$ etc.
 \end{proof}

 \begin{rk}\po
It is possible to define an almost complex structure by the above
approach on the bundles $T^{(r,s)}M=(TM)^{\ot r}\ot(T^*M)^{\ot
s}$, $S^{2i+1}TM$, $\La^{2j+1}TM$ etc (in some different manners),
but it won't be canonical (will depend on $\na$).
 \end{rk}

For two almost complex manifolds $(L,J_L)$ and $(M,J_M)$ a
canonical almost complex structure $\hat J$ on the space of
PH-1-jets
 $$
J^1_{PH}(L,M)=\{(x,y,\Phi)\,|\,x\in L, y\in M, \Phi\in T_x^*M\ot
T_yM:\ \Phi J_L=J_M\Phi\}
 $$
was introduced in \cite{Gau}. In particular, we get almost complex
structures on $J^1_{PH}(\C,M)=\C\times TM$ and
$J^1_{PH}(M,\C)=T^*M\times\C$. They are $\C$-translations
invariant and thus yield almost complex structures on $TM$ and
$T^*M$. Also the restriction of $\hat J$ defines a canonical
almost complex structure on $T^{(1,1)}_\C
M=\pi_{1,0}^{-1}(\Delta(M))$, where $\Delta(M)\subset M\times
M=J^0(M,M)$ is the diagonal PH-submanifold ($\pi_{1,0}:J^1_{PH}\to
J^0$ is the canonical projection). It can be shown that the
derived structures on $TM$ and $T^*M$ coincide with the ones
introduced above.

On the other hand, $J^1_{PH}(L,M)\subset T^{(1,1)}_\C(L\times M)$
is a PH-submanifold, whence the canonical structure $\hat J$ is a
generalization of that one from \cite{Gau}.

Note however that the higher PH-jet spaces $J^k_{PH}(L,M)$, $k>1$,
bear no structure in general (\cite{K1}) and usually are even
non-smooth.

 \vspace{4pt}
 {\bf Another canonical almost complex structure.}
An interesting issue is the paper \cite{LS}. An almost complex
structure on $TM$, which we denote $\check J$, is constructed
there via the deformation theory approach. It is not however new,
for it was introduced long before in \cite{YK} via the complete
lift operation $J\rightsquigarrow J^c$ (this fact was not noticed
in \cite{LS}). To see the coincidence $\check J=J^c$, note that in
local coordinates $(x^i,y^i)$ on $TM$ both structures have the
form $\begin{pmatrix}J&0\\\p J&J\end{pmatrix}$, where $\p=\sum
y^i\p_{x^i}$. It follows from \cite{YK,YI} that the structure
$\check J$ enjoys the same properties as the structure $\hat J$ in
theorem \ref{thm1}.

The structures $\hat J$ and $\check J$ differ because if we let
$Z$ denote multiplication by the complex number $z=a+ib\in\C$
along the fibers of $TM$, $Z(x,y)=(x,ay+bJy)$, we get $[\check
J,Z_*]=\!\!\begin{pmatrix}0&0\\bN_J(\p,\cdot)&0\end{pmatrix}$,
while from the very construction $\hat J\circ Z_*=Z_*\circ\hat J$.
Thus $\hat J\ne\check J$ unless $J$ is integrable.

We can also obtain $\hat J\not\equiv J^c$ from \cite{YI}, where
they provide a construction of almost complex structure $J^{\small
H}$ on $TM$ via horizontal lift of the connection $\bar
\nabla_XY=\nabla_YX+[X,Y]$ (equivalently
$\bar\Gamma^k_{ij}=\Gamma^k_{ji}$ in terms of Christoffel
symbols). By the construction $\hat J=J^{\small H}$ iff
$\bar\nabla$ is minimal and by the results in \S2.4 of \cite{YI}
$J^c=J^{\small H}$ iff $\nabla J=0$ (beware, without one of these
specifications the horizontal lift $J^{\small H}$ is
connection-dependent). But if $\nabla$ is an almost complex
connection, then $N_J(X,\cdot)=\bar\nabla_{JX}J-J\bar\nabla_XJ$,
whence $J^c=J^{\small H}\ne\hat J$ unless $N_J=0$.

The argumentation in \cite{LS} that $[\check J,Z_*]\ne0$ is
indirect and based on the fact that kernel of the Gromov operator
$D_u$ is not $J$-invariant. In \S\ref{sec6} we describe this
operator in terms of a canonical almost complex structure $\check
J$ on the normal bundle to a PH-curve (in fact, as notation
suggests, there is a relation between introduced canonical
structures on tangent and normal bundles).

 \begin{rk}\po
In \cite{YI} various lifts to tangent and cotangent bundles are
discussed. The complete lift of $J$ to the cotangent bundle is not
almost complex, but this is amended \cite{Sa} via the calibration
$J^c-\frac12\gamma(JN_J)$. The transformation is surprisingly
similar to our formula (\ref{e:J}) below, though we observe no
precise relations.
 \end{rk}

Let us call TB-I and TB-II the (total space of) tangent bundle
$TM$ equipped with the almost complex structures $\hat J$ or
$\check J$ respectively.

\section{\hps Almost complex normal bundle}
\label{sec2}

 \hspace{13.5pt}
Topologically the normal bundle $N_LM$ of a submanifold $L\subset
M$ is defined by the exact sequence:
 \begin{equation}\label{defNB}
0 \to TL \to TM\vert_L \to N_L M \to 0.
 \end{equation}

If $L$ is a complex submanifold of a complex manifold $M$, then
$N_LM$ is a holomorphic vector bundle over $L$ (the total space
and the projection are holomorphic, as well as fiberwise addition
and multiplication by complex numbers). In almost complex case
this is no longer so.

Let $\pi:N_LM\to L$ denote the projection and $\rho:L\to N_LM$ the
zero section.

 \vspace{4pt}
 {\bf NB-I structure.}
Here we apply the construction of \S\ref{sec1} to get a canonical
almost complex structure $\hat J$ on $N_LM$, called NB-I in what
follows:

 \begin{theorem}\po\label{thm2}
There exists a canonical almost complex structure $\J$ on the
total space of the normal bundle $N_LM$ to a PH-submanifold
$L\subset M$ such that:
 \begin{enumerate}
  \item
The maps $\pi:N_LM\to L$ and $\rho:L\to N_LM$ are
pseudoholomorphic.
  \item
The structure $\hat J$ is integrable iff $J|_L$ is integrable and
the $J$-antilinear by each argument part of the curvature
vanishes, $R_\na^{--}(X,Y)=0$, $\forall X,Y\in TL$, for some
minimal connection $\na$ totally geodesic and flat on $L$.
   \end{enumerate}
 \end{theorem}

 \begin{rk}\po
If $J|_L$ is integrable, the specified connection always exists
locally (the above integrability criterion is indeed local) and
then $R_\na^{--}(X,Y)$ does not depend on its choice (see appendix
\ref{a_A}). Moreover, $R_\na^{--}(X,Y)=0$ $\forall X,Y\in TL$,
whenever $J$ is integrable along $L$ to the second order:
$N_J(x)=0$ $\forall x\in L$.
 \end{rk}

 \begin{proof}
Let $\na$ be a minimal connection on $M$. It can be chosen so that
$L$ is totally geodesic. In fact, one chooses any linear
connection for which parallel transports along $L$ preserve $TL$
and note that the procedures of making the connection almost
complex and then minimal (see appendix \ref{a_A}) do not destroy
the property of $L$ to be totally geodesic.

We define a connection $\hat\nabla$ on the bundle $N_LM$ via
parallel transports as follows. Let $v=[\te]\in (N_LM)_x$ be the
class of $\te\in T_xM$ and let $\g(t)\subset L$ be a curve,
$\g(0)=x$. Calculate the parallel transport $\te(t)$ of $\te$
along $\g(t)$. Then define $v(t)=[\te(t)]$ to be the parallel
transport of $v$ along $\g(t)$. Since $L$ is totally geodesic, the
definition is correct ($\hat\nabla$-parallel transport of 0 is 0).
Moreover the connection $\hat\nabla$ is $\R$-linear. So as usual
in the theory of generalized connections we conclude that
$\hat\nabla$ is a linear connection.

Let $T_a(N_LM)=H_a\oplus V_a$ be the splitting into the horizontal
and vertical components induced by $\hat\nabla$, $a\in N_LM$. The
first space $H_a\stackrel{\pi_*}\simeq T_xL$ has a canonical
complex structure $J_1$ induced from $J|_L$ by $\pi_*$,
$x=\pi(a)$, and the second $V_a\simeq T_xM/T_xL$ inherits a
canonical complex structure $J_2$ from $J$ as the quotient. So we
obtain the structure $\hat J=J_1\oplus J_2$ on $T_a(N_LM)$ for
each $a$.

The same arguments as in Theorem \ref{thm1} show that the almost
complex structure $\hat J$ on $N_LM$ does not depend on the choice
of a minimal connection $\nabla$, preserving $TL$. The first
property of $\hat J$ is obvious. For the other one we use
 \begin{lem}\po
If a vector $Y\in T_a(N_LM)$ is vertical, then $N_{\hat
J}(\cdot,Y)=0$.
 \end{lem}
Actually, the fiber is integrable, so it is enough to consider the
pairing $N_{\hat J}(\hat X,Y)$, where $\hat X$ is the
$\hat\na$-lift of $X\in TL$. Recall (\cite{KN}) that $\hat\na_XY$
coincides with the Lie derivative $L_{\hat X}\ti Y$ of the section
$Y$ extended by translations to a vertical vector field $\ti Y$ on
$N_LM$ ($\hat X$ is the $\hat\nabla$-lift of any vector field
extending $X$; the result will not depend on an extension). Thus
$\hat\na_XY=[{\hat X},\ti Y]$ and we have (see also the remark
after proposition \ref{toj}):
 $$
N_{\hat J}(\hat X,Y)=\hat\na_{\hat JX}\hat JY-\hat J\hat\na_{\hat
JX}Y-\hat J\hat\na_X\hat JY -\hat\na_XY= (\hat\na_{\hat JX}\hat
J)Y+(\hat\na_{X}\hat J)\hat JY=0.
 $$

Now since the curvature of $\hat\na$ is
$R_{\hat\na}(X,Y)a=\widehat{[X,Y]}_a-[\hat X,\hat Y]_a$, we get:
 \begin{equation}\label{na}
N_{\hat J}(\hat X,\hat Y)_a=\widehat{N_J(X,Y)}_a
+4R_{\hat\na}^{--}(X,Y)a,\qquad X,Y\in TL.
 \end{equation}
For an integrable $J|_L$ we can choose minimal $\na$ to be flat on
$L$ and preserving $TL$, whence we get
$R_{\hat\na}^{--}(X,Y)=R_\na^{--}(X,Y)$ and the claim follows.
 \end{proof}

 \vspace{4pt}
 {\bf NB-II structure.}
From the integrability condition of Theorem \ref{thm2} we read off
that some features of 1-jet of the almost complex structure $J$
along $L$ are lost in $\hat J$ on $N_LM$. It is however possible
to keep most of them with another definition of the normal bundle
structure $\check{J}$, which we call NB-II.

 \begin{theorem}\po\label{thrm5}
There exists a canonical almost complex structure $\check{J}$ on
the total space of the normal bundle $N_LM$ to a PH-submanifold
$L\subset M$ such that:
 \begin{enumerate}
  \item
The maps $\pi:N_LM\to L$ and $\rho:L\to N_LM$ are
pseudoholomorphic.
 \item The structure $\check{J}$ is integrable iff the following 3
conditions hold:
   \begin{itemize}
   \item[$-$]$J|_L$ is integrable,
   \item[$-$]$(M,J)$ is normally integrable along $L$, i.\,e.\
$N_J(TL,TM|_L)\subset TL$,
   \item[$-$]The normal component $N_J^\perp=\chi\circ N_J$ vanishes
on $TL$ to the second order, where $\chi:TM|_L\to N_LM$ is the
natural projection.
   \end{itemize}
 \end{enumerate}
 \end{theorem}

 \begin{proof}
We describe the structure $\check J$ on the germ of zero section
in $N_LM$, which then uniquely determines it on the whole total
space.
Let $\ok^M_L$ be a tubular neighborhood of $L\subset M$. Fix a
$J$-invariant subbundle $F\subset TM|_L$ such that $TL\oplus
F=TM|_L$ (the totality of all such subspaces $F$ forms a bundle
over $L$ with contractible fibers). We identify $F=TM|_L/TL\simeq
N_LM$.

Let us fix some minimal connection $\nabla$ on $M$ with $L$ being
totally geodesic. Denote by
$N_LM\supset\ok_L^N\stackrel{\vp}\to\ok_L^M$ the
$\nabla$-exponential map that associates to the vector $v\in F_x$,
$x\in L$, the value $\gamma(1)$ along the $\nabla$-geodesic
$\gamma$ with initial conditions
$(\gamma(0),\dot\gamma(0))=(x,v)$.

Denote by $R^t$ the $t$-times dilatation $v\mapsto tv$ along the
fibers of $F$. We define:
 \begin{equation}\label{complPH}
 J^\vp=\vp_*^{-1} J \vp_*,\quad
 J_t=\op{ad}_{R^t}(J^\vp)=R^{1/t}_*J^\vp R_*^t\quad
 \text{ and }\quad \check{J}=\lim_{t\to0}J_t.
  \end{equation}
Consider local split coordinates $(x,y)$ on $N_LM$ such that
$L=\{y=0\}$ and the fibers of $F$ equal $\{x=\op{const}\}$. In
terms of these coordinates the limit process transforms the matrix
of $J^\vp$ as follows:
 $$
J^\vp=\begin{pmatrix}A(x,y) & C(x,y) \\ B(x,y) & D(x,y)
\end{pmatrix} \mapsto \check J=
\begin{pmatrix}A(x,0) & 0 \\ d_FB(x,y) & D(x,0),
\end{pmatrix}
 $$
where $d_FB(x,y)=\lim\limits_{t\to0}B(x,ty)/t$ (notice that
$B(x,0)=0$ because $TL$ is $J$-invariant).

Let us check independence of $\check J$ on $\nabla$ and $F$. When
we change the connection or the $J$-invariant subbundle, it is
equivalent to changing the map $\vp$ to $\tilde\vp$. In the above
split coordinates $(x^i,y^j)$ on $N_LM$ we have (assuming the
standard rule of summation by repeated indices)
 $$
\vp^{-1}\tilde\vp:(x^i,y^j)\mapsto (x^i+\a^i_k(x)y^k,y^j)+o(|y|)
 $$
(choice of the norm in $o(|y|)$ is not essential). Thus writing
the matrix of $J^\vp$ in block form we observe that the
transformation $J^\vp\mapsto J^{\tilde\vp}$ has the following
matrix form:
 \begin{equation}\label{block}
\begin{pmatrix}A & C \\ B & D \end{pmatrix}
\mapsto \Delta^{-1}\cdot
\begin{pmatrix}A & C \\ B & D \end{pmatrix}
\cdot\Delta= \begin{pmatrix}\tilde A & \tilde C
\\ \tilde B & \tilde D \end{pmatrix},
 \end{equation}
 where
 \begin{equation*}
 \Delta=d\bigl(\vp^{-1}\tilde\vp\bigr)=\1+
\begin{pmatrix}U & V \\ 0 & W \end{pmatrix}+o(|y|),
 \end{equation*}
and $U,W=o(1)$ have to vanish on $L$, but $V$ needs not to.

Since $B(x,0)=0$ we deduce from (\ref{block}): $\tilde
A(x,0)=A(x,0)$ and $\tilde D(x,0)=D(x,0)$. The transformation of
$C$ is inessential and $B$ changes to $\tilde
B(x,y)=B(x,y)(\1+o(1))$. Thus $d_F\tilde B(x,y)=d_FB(x,y)$ and we
see that the limit process (\ref{complPH}) gives a well-defined
result.

In addition we observe that the structure $\check J$ has affine
behavior w.r.t. $y$ and thus its restriction to $\ok_L^N$
determines the structure on the whole $N_LM$.

To prove integrability criterion we note that
$N_{\check{J}}=\lim\limits_{t\to\infty}N_{J_t}=
\lim\limits_{t\to\infty}\op{ad}_{R^t}(N_J)$. Consider $(x^i,y^j)$
as coordinates on both $\ok_L^N$ and $\ok_L^M$ using the
identification $\vp$.

Denote by $N_J^\perp$ the $y$-component of the value of $N_J$.
Note that $N_J^\perp$ is well-defined along $L$ and whenever
$J|_L$ is integrable, i.e. $N_J|_{TL}\equiv0$, its 1-jet is
well-defined. Then we calculate:
 \begin{equation}\label{apa}
N_{\check{J}}(\p_{x^i},\p_{x^j})=N_J(\p_{x^i},\p_{x^j})|_{y=0}+
y^k\bigl(\p_{y^k}N_J^\perp(\p_{x^i},\p_{x^j})|_{y=0}\bigr)
 \end{equation}
and
 \begin{equation}\label{ape}
N_{\check{J}}(\p_{x^i},\p_{y^j})=N^\perp_J(\p_{x^i},\p_{y^j})|_{y=0},\quad
N_{\check{J}}(\p_{y^i},\p_{y^j})=0,
 \end{equation}
The claim follows.
 \end{proof}

 \vspace{4pt}
If $\op{codim}_\C L=1$, then the connection $\nabla$ can be chosen
so that the exponential image of the vertical foliation $\vp(F)$
is $J$-holomorphic. This follows from

 \begin{prop}\po\label{prop:1}
Small neighborhood $\ok_L$ of a PH-submanifold $L^{2n-2}\subset
M^{2n}$ 
can be foliated by transversal PH-disks $D^2$.
 \end{prop}

 \begin{proof}
This follows from Nijenhuis-Woolf theorem~\cite{NW} on the
existence of a small PH-disk in a given direction, smoothly
depending on it.
 \end{proof}

 \begin{rk}\po
For $n=2$ a construction of certain structure $\bar J$ on
$\ok_L^M$, using the dilatation $R^t$ and based on the idea of
Proposition \ref{prop:1}, was used in \cite{M2}.
 \end{rk}

Denote by $N_L^IM$ and $N_L^{II}M$ the normal bundle equipped with
the NB-I structure $\hat J$ or with the NB-II structure $\check J$
respectively. The tangent bundle structures TB-I and TB-II can be
deduced from the normal ones via the diagonal embedding
$\Delta:M\hookrightarrow M\times M$ because $N_{\Delta(M)}(M\times
M)\simeq TM$.

We are going to relate the concept of NB-II with the deformation
theory. The following statement will be used in \S\ref{sec6}.

 \begin{prop}\po\label{proo}
Let $\phi_t:(\CC,J_\CC^t)\to(M,J_M)$ with $\phi_0(\CC)\subset L$
be a family of $J$-holomorphic embeddings. Then
$\phi_t'|_{t=0}:(\CC,J_\CC^0)\hookrightarrow N_L^{II}M$ is a
PH-embedding. In particular, deformations of $\CC=L$ lead to
PH-sections $\phi_t'|_{t=0}:(L,J_L^0)\hookrightarrow N_L^{II}M$.
 \end{prop}

Notice that by virtue of the relation between NB-I and NB-II from
the next section the embedding $\phi_t'|_{t=0}$ of $L$ into
$N_L^IM$ is not pseudoholomorphic.

 \begin{proof}
We have $J_M\, d\phi_t=d\phi_t J^t_\CC$, whence
 $$
(R^{1/t}_* J_M R^t_*) (R^{1/t}_* d\phi_t)=(R^{1/t}_* d\phi_t)
J^t_\CC.
 $$
In the limit $t\to0$ we get: $\check J\, d\phi'_0=d\phi'_0\,
J^0_\CC$.
 \end{proof}

 \vspace{4pt}
This proposition leads to an equivalent definition of the NB-II
structure $\check J$.

Consider $x\in L$, $v\in(N_LM)_x$ and $\zeta\in T_v(N_LM)$. Let
$w\in T_xM$ represent $v$, $v=[w]$. Consider a curve $\g(t)$ in
$M$ with $\g(0)=x$, $\dot\g(0)=w$ and a vector field along the
curve $\xi_t\in T_{\g(t)}M$ that represents $\zeta$. Then
$\eta_t=J_M\xi_t\in T_{\g(t)}M$ represents $\varsigma=\check
J\zeta\in T_v(N_LM)$.

In fact, there exists a family of PH-disks
$\phi_t:(D^2_\ve,J_0)\to(M,J_M)$ with $\phi_t(0)=\g(t)$,
$d_0\phi_t(1)=\xi_t$. Then $d_0\phi_t(i)=\eta_t\in T_{\g(t)}M$,
where $1,i\in T_0D^2_\ve$.

From this alternative definition we obtain

 \begin{lem}\po\label{prop8}
Let $L\subset M$ be a PH-submanifold w.r.t. two almost complex
structures  $J_1$ and $J_2$ with equal normal bundles $N_L^{II}M$.
Then $\nabla_Y(J_1-J_2)(X)=0$ for all $X\in TL$ and $Y\in TM|_L$
(the choice of connection is inessential).
 \end{lem}

 \begin{proof}
Let $\g(t)\subset M$ be a curve with $\g(0)=x\in L$,
$\dot\g(0)=Y$. Consider two family of PH-disks
$\phi_t:(D^2_\ve,J_0)\to(M,J_1)$ and
$\psi_t:(D^2_\ve,J_0)\to(M,J_2)$ with $\phi_t(0)=\psi_t(0)=\g(t)$
and $d_0\phi_0(1)=d_0\psi_0(1)=X$. We can suppose that they induce
the same map $\phi'_0=\psi'_0:(D^2_\ve,J_0)\to N_L^{II}M$. Then:
 $$
\nabla_Y(J_1-J_2)(X)= \left.\tfrac
d{dt}\right|_{t=0}\bigl(J_1d_0\phi_t(1)-J_2d_0\psi_t(1)\bigr)
=d_0\phi'_0(i)-d_0\psi'_0(i)=0.
 $$
 \end{proof}

\section{\hps Pseudoholomorphic vector bundles}
\label{sec3}

 \hspace{13.5pt}
Consider a real vector bundle $\pi:(E,\hat J)\stackrel{F}\to(L,J)$
with almost complex total space, base and projection: $\pi_*\hat
J=J\pi_*$. The following statement is obvious:

 \begin{prop}\po\label{lem:3}
The Nijenhuis tensor $N_{\hat J}$ is projectible: $\pi_*N_{\hat
J}=N_J\circ\La^2\pi_*$.\qed
 \end{prop}

 \begin{cor}\po\label{999}
Let $(L,J)$ be integrable (for example $\dim_\C L =1$). Then we
have: $\op{Im}(N_{\hat J})\subset TF$.\!\qed
 \end{cor}

 \begin{dfn}\po\label{fed}%
Call $\pi$ a {\em almost holomorphic vector bundle\/} (we write PH
-- pseudoholomorphic), if the restrictions $\hat J|_{F_x}$ are
constant coefficients complex structures on the fibers and there
exists a linear (not necessary $J$-linear) connection $\hat\nabla$
on $\pi$ such that the $\hat\nabla$-lift $\C$-splits the exact
PH-sequences
 $$
0\to F_x\to T_aE\stackrel\dashleftarrow\longrightarrow
T_xL\to 0, \qquad x=\pi(a),
 $$
In this case the zero section $L\subset E$ is a $\hat
J$-holomorphic submanifold.
 \end{dfn}

 \begin{prop}\po
The canonical almost complex structures $\hat J,\check J$ on $TM$
and $\hat J,\check J$ on $N_LM$ are PH vector bundle structures.
 \end{prop}

 \begin{proof}
For TB-I and NB-I structures $\hat J$ the claim is implied
directly by the construction. For TB-II and NB-II structures
$\check J$ this follows from the explicit formulas and the affine
behavior by the fiber coordinates.
 \end{proof}

 \vspace{4pt}
Consider an arbitrary splitting $TE=H\oplus V$ into horizontal and
vertical components. Restricting the first argument of the
Nijenhuis tensor to $H$ and the second to $V=TF$ we obtain a
tensor $N'{\!\!}_{\!\hat J}: \pi^*TL\ot TF\to TF$.

 \begin{prop}\po\label{135}
The tensor $N'{\!\!}_{\!\hat J}$ does not depend on a choice of
horizontal component $H$ (not necessary $\hat J$-lift) and is
constant along the fibers. So it is lifted from a canonical tensor
(we will use the same notation) $N'{\!\!}_{\!\hat J}:TL\ot F\to F$
with $\hat J$-invariant image $\Pi'{\!}_{\!\hat J}=N_{\hat
J}(H,V)\subset F$.
 \end{prop}

 \begin{proof}
Independence of $H$ follows from Proposition \ref{lem:3}. Let us
prove constancy along the fibers $F$. Let $\hat\na$ be a
connection from the definition.

Denote $\# j=j-(-1)^j$. There are local coordinates $(x^i,y^j)$ on
$\pi^{-1}(U)=U\times F$, with $x$ a base coordinate and $y$ a
linear fiber coordinate, such that the structure ${\hat J}|_F$ has
constant coefficients w.r.t.\ $y$:
 \begin{equation}\label{nyn}
{\hat J}\p_{y^j}=(-1)^{j-1}\p_{y^{\# j}}.
 \end{equation}
Let $\hat\na_{\p_{x^i}}\p_{y^j}=\G_{ij}^k(x)\p_{y^k}$. The
$\hat\na$-lift of $\p_{x^i}$ is:
$\hat\p_{x^i}=\p_{x^i}-\G_{ij}^s(x)y^j\p_{y^s}$.

Let $J\p_{x^i}=a_i^k(x)\p_{x^k}$ on the base. Then $\hat
J\hat\p_{x^i}=a_i^k(x)\hat\p_{x^k}$ and we get:
 \begin{equation}\label{nan}
\hat J\p_{x^i}=a_i^k\p_{x^k}+
\bigl((-1)^s\G_{ij}^{\#s}-a_i^k\G_{kj}^s\bigr)y^j\p_{y^s}.
 \end{equation}
Thus $N_{\hat J}(\p_{x^i},\p_{y^j})=\g_{ij}^s(x)\p_{y^s}$ is
expressed via the Christoffel symbols as
 \begin{equation}\label{nana}
\g_{ij}^s=(-1)^{s+j}\G_{i,\#j}^{\#s}-(-1)^sa_i^k\G_{kj}^{\#s}
-(-1)^ja_i^k\G_{k,\#j}^s-\G_{ij}^s,
 \end{equation}
so it is constant along the fibers. Note that
$\op{rk}(\Pi'{\!}_{\!\hat J})$ can vary with $x\in L$.
 \end{proof}

 \begin{dfn}\po\label{touj}
Let us call a PH-bundle almost complex structure $\hat J$ on
$(E,\pi)$ {\em normally integrable\/} if $N'{\!\!}_{\!\hat J}=0$.
 \end{dfn}

For such a structure integrability is equivalent to integrability
of $(L,J)$ and vanishing of $R_{\hat\na}^{--}$ (cf. proofs of
Theorems \ref{thm2}, \ref{thrm5} and formula (\ref{na})). In
particular:

 \begin{prop}\po\label{toj}
Normally integrable PH bundles over holomorphic curves are
holomorphic. \qed
 \end{prop}

If $\hat\na$ is obtained from a minimal connection $\na$, as for
the structures $\hat J$ of \S\ref{sec1}-\ref{sec2}, then it
additionally preserves $\hat J|_F$, meaning
$\G_{i,\#j}^{\#s}=(-1)^{s+j}\G_{ij}^s$. So (\ref{nana}) implies
$\gamma_{ij}^s=0$ and $N'{\!\!}_{\!\hat J}=0$. In particular, the
NB-I structures $\hat J$ over a PH-curve is normally integrable
(while the NB-II structure $\check J$ is usually not). To describe
such structures in general notice that formula (\ref{na}) implies
the following:

 \begin{prop}\po\label{246}
If a PH bundle structure $\hat J$ is normally integrable, then
restriction of the Nijenhuis tensor to both horizontal components
determines a canonical tensor $N''{\!\!}_{\!\hat
J}:\pi^*\La^2TL\to TE$ with the image $\Pi''{\!}_{\!\hat J}
=N_{\hat J}(H,H)\subset TE$ being a $\hat J$-invariant
differential system. This tensor projects to the tensor $N_J$ on
the base and is affine-linear along the fiber. \qed
 \end{prop}

Let $a\in E$ and $x=\pi(a)\in L$ be its projection. Denote by
$r=r_a\in F_x\subset T_aE$ the radius-vector $\vec{xa}$.

 \begin{theorem}\po\label{th:1}
Let $(E,\hat J,\pi)$ be a pseudoholomorphic vector bundle over an
almost complex manifold $(L,J)$. Then $\hat J$ can be expressed
via some normally integrable PH vector bundle structure $J_0$ and
the tensor $N_{\hat J}$ by the formula:
 \begin{equation}\label{222}
\hat J=J_0+\frac12J_0 N_{\hat J}(r,\cdot).
 \end{equation}
 \end{theorem}

 \begin{proof}
Let us define the structure by the formula
 \begin{equation}
J_0=\hat J-\frac12{\hat J}N_{\hat J}(r,\cdot).
 \label{e:J}
 \end{equation}
Since $N_{\hat J}|_F\equiv0$ this structure $J_0|_F=\hat J|_F$ is
a constant complex structure on the fibers $F$, proving formula
(\ref{222}) for $\hat J$.

To show that the structure $J_0$ is almost complex, we note that
$N_{\hat J}(r,Y)\in F$ for any $Y$ and $N_{\hat J}(r,Y)=0$ for
$Y\in F$. Therefore
 $$
J_0^2={\hat J}^2- \frac12{\hat J}^2 N_{\hat J}
(r,\cdot)-\frac12{\hat J}N_{\hat J}(r,\hat J\cdot)+ \frac14{\hat
J}N_{\hat J}(r,{\hat J}N_{\hat J}(r,\cdot))={\hat J}^2=-\1.
 $$

To obtain $N'{\!\!}_{J_0}=0$ we use (\ref{e:J}) and the
coordinates of proposition \ref{135}:
 $$
 \begin{array}{rl}
&\hspace{-15pt} N_{J_0}(\p_{x^i},\p_{y^j})=\\
\vphantom{\left\{\dfrac{1^1}{2^2}\right\}}
&\hspace{-10pt}= N_{\hat
J}(\p_{x^i},\p_{y^j})-[\frac12{\hat J}N_{\hat J}
(y^s\p_{y^s},\p_{x^i}),(-1)^{j-1}\p_{y^{\# j}}]
+\hat J[\frac12{\hat J}N_{\hat J}(y^s\p_{y^s},\p_{x^i}),\p_{y^j}]\\
&\hspace{-10pt}=N_{\hat J}(\p_{x^i},\p_{y^j}) +\frac12(-1)^j\hat J
N_{\hat J}(\p_{x^i},\p_{y^{\# j}}) +\frac12\hat J^2 N_{\hat
J}(\p_{x^i},\p_{y^j})=0,
 \end{array}
 $$
where we expressed $r=y^s\p_{y^s}$. The claim follows.
 \end{proof}

 \begin{cor}\po\label{coro}
If the base is a PH-curve, $\dim_\C L=1$, then the structure $J_0$
in formula (\ref{222}) is complex analytic, making $\pi$ into a
holomorphic vector bundle. \qed
 \end{cor}

 \begin{dfn}\po
Let us call the structure $J_0$ of theorem \ref{th:1} the normally
integrable form {\em (n.i.f.)\/} of the PH-bundle structure $\hat
J$.
 \end{dfn}

Certainly normally integrable form of a normally integrable
structure (e.g. TB-I or NB-I) $\hat J$ is this structure itself.
Now we will describe a relation between NB-I and NB-II structures
(implying a similar relation for TB-I and TB-II). We consider the
latter as a general pseudoholomorphic vector bundle.

 \begin{theorem}\po\label{t14}
Let $(L,J)$ be the zero section of a PH vector bundle $(E,\hat
J,\pi)$. Then its NB-I structure coincides with the n.i.f.\ $J_0$
of the structure $\hat J$ as in (\ref{222}).
 \end{theorem}

 \begin{proof}
We use formulae (\ref{nyn}) and (\ref{nan}) for the almost complex
structure. Consider a linear connection $\na$, given by the
relations
 $$
\na_{\p_{x^i}}\p_{x^j}=0,\
\na_{\p_{x^i}}\p_{y^j}=\G_{ij}^k(x)\p_{y^k},\
\na_{\p_{y^i}}\p_{x^j}=0,\ \na_{\p_{y^i}} \p_{y^j}=0.
 $$
Calculate by it a minimal connection $\ti\na$ by the algorithm of
appendix \ref{a_A}. It in turn produces the following connection
on the normal bundle $N_LE\simeq E$:
 $$
\bar\na_{\p_{x^i}}\p_{y^j}=\Bigl(\frac38\G_{ij}^s
+\frac18(-1)^sa_i^k\G_{kj}^{\#s}-\frac18(-1)^ja_i^k\G_{k,\#j}^s
+\frac38(-1)^{s+j}\G_{i,\#j}^{\#s}\Bigr)\p_{y^s}.
 $$
Using the relation $\bar J\bar\p_{x^i}=a_i^k\bar\p_{x^k}$ we get
the formula
 \begin{equation*}
\bar J\p_{x^i}=a_i^k\p_{x^k}+\frac12\Bigl(
(-1)^s\G_{ij}^{\#s}-(-1)^j\G_{i,\#j}^s-a_i^k\G_{kj}^s
-(-1)^{s+j}a_i^k\G_{k,\#j}^{\#s}\Bigr)y^j\p_{y^s},
 \end{equation*}
which together with the formula $\bar J|_F=\hat J|_F$ (\ref{nyn})
describes the NB-I structure $(E,\bar J)$ of the zero section.

But substitution of formulae (\ref{nan}) and (\ref{nana}) into
(\ref{e:J}) yields the same expressions for $J_0$, proving the
claim: $\bar J=J_0$.
 \end{proof}

 \vspace{4pt}
Thus the two PH-bundles $N_L^IM$ and $N_L^{II}M$ are related as
follows:

 \begin{equation}\label{65}
\text{NB-II}\ \stackrel{\text{n.i.f.}}\longrightarrow\ \text{NB-I}
 \end{equation}

 {\bf Relation to other generalizations of holomorphic bundles.}
Our PH-vector bundle structures differ from "bundle almost complex
structures" of Bartolomeis and Tian \cite{BT}, because (see
\S\ref{sec1}) the multiplication morphism $\mu:\C\times E\to E$ is
not pseudoholomorphic in general (though its restriction
$\mu:\R\times E\to E$ is). But they satisfy the requirements of
"almost holomorphic vector bundles" by Lempert and Sz\"oke
\cite{LS}. Actually our definitions are equivalent:

 \begin{prop}\po
$(E,M,\pi)$ is a PH vector bundle structure iff the
fiber-wise addition $\a: E\times_M E\to E$ is a PH-map.
 \end{prop}

 \begin{proof}
The almost complex structure $\hat J$ on $E\times_M E$ is induced
from the natural product structure on $E\times E$, since the
former is the preimage of the diagonal $\Delta(M)\subset M\times
M$ (which is pseudoholomorphic).

In local coordinates $(x^i,y^j)$ the structure $\hat J$ on $E$ is
given by formulae (\ref{nyn})-(\ref{nan}). Then the structure on
$E\times_M E$ is given in local coordinates $(x^i,z^j,w^k)$ as
follows (we do not specify coefficients $b_{ij}^s$ via the
Christoffel coefficients):
 \begin{gather*}
\hat
J\p_{x^i}=a_i^k(x)\p_{x^k}+b_{ij}^s(x)z^j\p_{z^s}+b_{ij}^s(x)w^j\p_{w^s},\\
{\hat J}\p_{z^j}=(-1)^{j-1}\p_{z^{\# j}},\quad {\hat
J}\p_{w^j}=(-1)^{j-1}\p_{w^{\# j}}.
 \end{gather*}
The map $\a_*$ maps both $\p_{z^j}$ and $\p_{w^j}$ to $\p_{y^j}$.
It is enough to check that it is a PH-map only on the basic
vectors. Consider a point $(x,z,w)\stackrel\a\mapsto(x,y=z+w)$.
For $\p_{z^j}$ and $\p_{w^j}$ we have: $\a_*\hat J=\hat J\,\a_*$.
And for the horizontal vectors:
 $$
\hat J_{(x,y)}\a_*(\p_{x^i})-\a_*(\hat J_{(x,z,w)}\p_{x^i})=
b_{ij}^s(x)(y^j\p_{y^s}-z^j\p_{y^s}-w^j\p_{y^s})=0.
 $$
Thus if $(E,\hat J)$ is a PH bundle, the map $\a$ is
pseudoholomorphic.

On the other hand if $\a$ is a PH-map, then the above arguments
show local existence of a connection $\hat\nabla$,
satisfying the requirement of definition \ref{defNB}. The space of
such connections is contractible, whence the global existence.
 \end{proof}

\section{\hps Normal form of 1-jet of $J$ along a submanifold}
\label{sec4}

 \hspace{13.5pt}
Consider the ideal of $\R$-valued functions corresponding to a
submanifold $L$:
 $$
\mu_L=\{f\in C^\infty(M)\,|\,f(L)=0\}.
 $$
Its degrees determine the filtration $\mu^k$ on every
$C^\infty(M)$-module, in particular we can talk about jets of
tensor fields along $L$:
$J^k(\mathcal{T})=C^\infty(\mathcal{T})/\mu_L^{k+1}
C^\infty(\mathcal{T})$.

 \begin{theorem}\po\label{th:2}
Let $L\subset M$ be a PH-submanifold with respect to two almost
complex structures $J_1$ and $J_2$. Assume that the following
holds:
 \begin{enumerate}
 \item For every point $x\in L$: $J_1(x)=J_2(x)$,
$N_{J_1}(x)=N_{J_2}(x)$.
 \item The normal bundles $N_L^{II}M$ w.r.t. the structures $J_1$ and $J_2$
coincide.
 \end{enumerate}
Then $J_1$ and $J_2$ are 1-jet equivalent along $L$: There exists
a diffeomorphism $\vp$ of a neighborhood $\ok(L)$, such that
$\vp|_L=\op{Id}$, $d_x\vp=\1$ for all $x\in L$ and
 $$
J_2=\vp^*J_1\,\op{mod}\mu_L^2.
 $$
 \end{theorem}

Notice that the required conditions are necessary for
1-equivalence.

 \begin{rk}\po\label{rem_ns}
When $J_i$ are integrable and defined on different manifolds
$M_i$, but with the same normal bundle $N$, there is the
Nirenberg-Spencer cohomology obstruction $\op{ns}_0(J_1,J_2)\in
H^1(L;TL\ot N^*)$ (\cite{NS,MR}) for the 1$^\text{st}$ order
equivalence. It equals the difference of obstructions to splitting
the normal bundle sequence (\ref{defNB}). In particular, if the
sequences are isomorphic, then $\op{ns}_0(J_1,J_2)=0$.

In our case $M_1=M_2$ and the class $\op{ns}_0$ vanishes by
condition 2. However if we want to formulate the equivalence of
1-jets of $J_1$ and $J_2$ on different manifolds, we should
require $\op{ns}_0(J_1,J_2)=0$, where the latter will be
determined via NB-I structure (common for $J_1$ and $J_2$) and
sequence (\ref{defNB}).
 \end{rk}

In the calculations below we denote by $\stackrel{\textbf{.}}=$
the equivalence modulo $\mu_L$ (equality of 0-jets) and by
$\stackrel{\textbf{..}}=$ the equivalence modulo $\mu_L^2$
(equality of 1-jets).
 \vspace{4pt}

 \begin{proof}
Let us choose a minimal connection $\nabla$ near $L$ with
$L$ being totally geodesic. We wish to find $\vp:\ok_L\to\ok_L$
with $d\vp\circ J_1\stackrel{\textbf{..}}=J_2\circ d\vp$. This
implies
 \begin{equation}\label{NJNJ}
d\vp\circ N_{J_1}\stackrel{\textbf{.}}=N_{J_2}\circ\La^2d\vp.
 \end{equation}
Thus the tensor $\nabla d\vp$ is symmetric along $L$. Indeed, we
have: $(\na_Xd\vp)(Y)=\na_{d\vp(X)}(d\vp(Y))-d\vp(\na_XY)$ and so
 \begin{multline*}
(\na_Xd\vp)(Y)-(\na_Yd\vp)(X)\\
=T_\na(d\vp(X),d\vp(Y))+[d\vp(X),d\vp(Y)] -d\vp\,
T_\na(X,Y)-d\vp[X,Y]\\
=\tfrac14(N_{J_2}\circ\La^2d\vp-d\vp\circ
N_{J_1})(X,Y)\stackrel{\textbf{.}}=0.
 \end{multline*}
Denote $\Phi^{(2)}=\na d\vp\in C^\infty(S^2T^*M\ot TM|_L)$. In
terms of this tensor, the condition $d\vp\circ
J_1\stackrel{\textbf{..}}=J_2\circ d\vp$ holds iff for all $X,Y\in
TM|_L$ we have:
 \begin{equation}\label{diff}
\P^{(2)}(X,J_1Y)+d\vp\,(\nabla_X J_1)(Y)= J_2\P^{(2)}(X,Y)+
(\nabla_{d\vp\,(X)}J_2)(d\vp\,(Y)).
 \end{equation}
Denote
 \begin{equation}\label{e:P'}
P(X,Y)=(\nabla_{d\vp\,(X)}J_2)(d\vp\,(Y))-d\vp\,(\nabla_X J_1)(Y).
 \end{equation}
This yields the followings property along $L$:
 $$
P(X,J_1Y)=-J_2\c P(X,Y),
 $$
which implies that $P(X,Y)=J_2B(X,Y)-B(X,J_1Y)$ for some
$(2,1)$-tensor $B$. Conditions (\ref{e:P'}) and (\ref{NJNJ}) yield
(with $J=J_1=J_2$ along $L$):
 $$
P(X,Y)-P(Y,X)=P(JX,JY)-P(JY,JX).
 $$

From this we obtain a solution (similarly to Theorem~1 of
\cite{K1})
 \begin{multline*}
\Phi^{(2)}(X,Y)=-\frac12[B(X,Y)+B(Y,X)]\\
+\dfrac{J}4[B(JX,Y)+B(JY,X)-B(X,JY)-B(Y,JX)]
 \end{multline*}
of the equation $P(X,Y)=\P^{(2)}(X,J_1 Y)-J_2\P^{(2)}(X,Y)$ and
hence of (\ref{diff}).

We want to construct a map with $d\vp\stackrel{\textbf{.}}=\1$.
This requirement, equation (\ref{e:P'}) and assumptions of the
theorem imply that $P(X,Y)=0$ along $L$ if $X\in TL$ or $Y\in TL$
(see lemma \ref{prop8}). Thus we can choose $B$ with the same
property and get $\Phi^{(2)}(X,Y)=0$ if at leat one of $X,Y$
belongs to $TL$.

Now we integrate the symbols $\Phi^{(2)}$ to get the 2-jet of
$\vp$ along $L$, using the Taylor-Maclaurin decomposition by the
normal coordinate $y$ along a complimentary to $TL$ $J$-invariant
subbundle $F$, as in the proof of theorem \ref{thrm5}.
 \end{proof}

 \vspace{4pt}
A combination of Theorems \ref{th:1} and \ref{th:2} yields normal
forms of 1-jets of almost complex structures $J$ along a
PH-submanifold $L\subset(M,J)$.

Let us choose a $J$-invariant subbundle $F$ transversal to $L$, as
in the proof of Theorem \ref{thrm5}.
Consider the radial vector field $r$, which equals $\vec{xa}$ at
the point $a\in F_x$, $x\in L$ (as in Theorem \ref{th:1}). Let $A$
be an automorphism of $TM|_{\ok_L}$, which equals
$A=\begin{pmatrix} 1/2&0 \\ 0&1/4
\end{pmatrix}$ along $L$ in the decomposition $TM|_L=TL\oplus F$.

 \begin{theorem}\po\label{th:3}
Let $L\subset (M,J)$ be a PH-submanifold and $N_J\in C^\infty(
\Lambda^2T^*M\ot_{\bar\C}TM|_L)$ be the field of Nijenhuis tensors
of $J$ along it. Then there exist a normally integrable almost
complex structure $J_0$ in a neighborhood $\ok_L\subset M$ and a
diffeomorphism $\vp$ of $\ok_L$ such that $J_0=J$ along $L$,
$d\vp=\1$ along $L$ and we have:
 \begin{equation}
\vp^*J=J_0+J_0N_J(r,A\cdot)\,\op{mod}\mu^2_L.
 \label{e:MU}
 \end{equation}
In particular, when $L$ is a PH-curve, the structure $J_0$ can be
chosen complex.
 \end{theorem}

 \begin{proof}
Define $J'=J-JN_J(r,A\cdot)$. This is an almost complex structure
$\op{mod}\mu^2_L$ (see \cite{K1} about such jets). In fact,
$J'\stackrel{\textbf{.}}=J$ and $AJ\stackrel{\textbf{.}}=JA$, so
that
 $$
J'^2\stackrel{\textbf{..}}= J^2-J^2N_J(r,A\cdot)-JN_J(r,AJ\cdot)
\stackrel{\textbf{..}}=J^2=-\1.
 $$
Notice that we get $J\stackrel{\textbf{..}}=J'+J'N_J(r,A\cdot)$.

Let $J_0=\check J'$ be the corresponding NB-II structure (it is
already a genuine almost complex structure). Then $\ti
J=J_0+J_0N_J(r,A\cdot)$ is an almost complex structure
$\op{mod}\mu^2_L$ and it has the same NB-II structure as the
structure $J$.

Now we want to check the second part of assumption 1 in Theorem
\ref{th:2} for the structures $J,\ti J$ (we obviously have $\ti
J\stackrel{\textbf{.}}=J$).

Let $X^\perp$ denote the $F$-component of $X\in TM|_L$. Then we
get $[X,N_J(r,Y)]\stackrel{\textbf{.}}=N_J(X^\perp,Y)$ (compare
with the proof of Theorem \ref{th:1}, where $r=y^i\p_{y^i}$ in
local coordinates). And so we calculate:
 \begin{align*}
N_{J'}(X,Y)&\stackrel{\textbf{.}}=N_J(X,Y)
-[JX,JN_J(r,AY)]-[JN_J(r,AX),JY]\\
&\hspace{59pt} +J[X,JN_J(r,AY)]+J[JN_J(r,AX),Y]\\
&\stackrel{\textbf{.}}=N_J(X,Y)
-JN_J(JX^\perp,AY)-JN_J(AX,JY^\perp)\\
&\hspace{59pt} -N_J(X^\perp,AY)-N_J(AX,Y^\perp)\\
&=N_J(X,Y) -2N_J(X^\perp,AY)-2N_J(AX,Y^\perp).
 \end{align*}

Thus if $X,Y\in TL$, then $N_{J'}(X,Y)=N_J(X,Y)$. If $X\in TL$,
$Y\in F$, then $N_{J'}(X,Y)=N_J(X,Y)-2N_J(AX,Y)=0$. And if $X,Y\in
F$, then $N_{J'}(X,Y)=N_J(X,Y)-2N_J(AX,Y)-2N_J(X,AY)=0$.

Therefore, $N_{J'}$ vanishes for vertical vectors and $J_0$ is
normally integrable. In particular, $J_0$ is the NB-I structure of
the structure $J'$, see (\ref{65}).

By a calculation, similar to the above one, we obtain along $L$:
 $$
N_{\ti J}(X,Y)=N_{J_0}(X,Y) +2N_J(X^\perp,AY)+2N_J(AX,Y^\perp).
 $$
Since $N_{J_0}(X,Y)=0$ if $X$ or $Y$ belongs to $F$ and
$N_{J_0}|_{TL}=N_J|_{TL}$, we conclude that $N_{\ti
J}(X,Y)=N_J(X,Y)$ for all $X,Y\in TM|_L$.

Thus from Theorem \ref{th:2} we get a local diffeomorphism $\vp$
identical up to the first order on $L$ and such that $\ti
J\stackrel{\textbf{..}}=\vp^*J$.
 \end{proof}

 \begin{rk}\po
When $L$ is a point, the structure $J_0$ can also be chosen
complex. Moreover in this case $A=1/4$ and formula (\ref{e:MU})
looks especially simple. We write it in local coordinates $(x^i)$
centered at the given point $x_0\in M$:
 $$
J_i^k=(-1)^k\d_i^{\# k}-(-1)^k\tfrac14N_{ij}^{\# k}(0)x^j+o(|x|).
 $$
A general way to obtain similar formulae for jets at a point is
related to the structural function (Weyl tensor) of the
corresponding geometric structure (\cite{KL}).
 \end{rk}

\section{\hps Four-dimensional case and Arnold's question}
\label{sec5}

 \hspace{13.5pt}
In this and next sections we consider the special case $\dim M=4$.
Proper PH-submanifolds are PH-curves $L^2\subset(M^4,J)$. So
$N_L^IM=(N_LM,\hat J)$ is a holomorphic line bundle, while
$N_L^{II}M=(N_LM,\check J)$ is a PH-line bundle.

{\it Nijenhuis\- tensor characteristic distribution\/}
$\Pi=\op{Im}(N_J)\subset TM^4$ (\cite{K3}) is $J$-invariant and
has rank 2 in the domain of non-integrability for $J$, $N_J\ne0$.

 \begin{prop}\po
At the points $x\in L$, where the Nijenhuis tensor characteristic
distribution $\Pi$ is transversal to $L$, the same happens to the
NB-II characteristic distribution $\check\Pi$. But $N_{\check
J}(x)=0$ at the points $x$, where $\Pi\subset TL$.
 \end{prop}

 \begin{proof}
This follows from formulae (\ref{ape}).
 \end{proof}

 \begin{cor}\po
If the Nijenhuis tensor characteristic distribution $\Pi$ is
tangent to $L$, then the NB-I and NB-II structures coincide and
are holomorphic. \qed
 \end{cor}

Holomorphic line bundles over a genus $g$ curve $L=\Sigma_g$ are
parameterized by $g$ complex parameters. Line bundles over
rational curves $L=\bar\C\simeq S^2$ are determined by the
topological type, i.e.\ by the self-intersection number $L\cdot L$
of the zero section. But for other curves the holomorphic and
differentiable types of holomorphic bundles are different.

A holomorphic line bundle over an elliptic curve
$L=\C/\Z^2(2\pi,\oo)\simeq T^2$ ($g=1$), $\oo\in\C\setminus\R$,
depends on one parameter $\l\in\C\setminus\{0\}$. If the zero
section has self-intersection number $p$, the bundle is: $E\to
T^2$, $(z,w)\mapsto z$, $J_0=i$, with
 \begin{equation}
E=\C^2/(z,w)\sim(z+2\pi,w)\sim(z+\oo,\l e^{-ipz}w).
 \label{e1}
 \end{equation}

The pair $(\oo,\l)$ can 
be chosen to satisfy: $|\oo|\ge2\pi$, $-\pi<|\op{Re}\oo|\le\pi$,
$\op{Im}\oo>0$, $e^{-\op{Im}\oo}<|\l|\le1$. The number $\oo$ is
defined by the restriction $J_0|_{T^2}$ and the number $\l$ is
defined by 1-jet of the structure $J_0$ on $T^2$.

A PH-line bundle $(N_LM,\check J)$ over a genus $g$ curve
$L=\Sigma^2_g$ is parametrized by $g$ complex parameters (for NB-I
structure $J_0=\hat J$), a cohomology class $\op{ns}_0\in
H^1(L;TL\ot N^*)$, see remark \ref{rem_ns}, and a smooth 1-form
$N_{\check J}\in \O^1(L;\op{aut}_{\bar\C}(N_LM))$.

Consider an elliptic curve $L=T^2$ in a complex surface
$(M^4,J_0)$ with the normal bundle $N_L^IM$ given by (\ref{e1}).
For $p<0$ (\cite{Gra}) or $p=0$ and generic pair $(\oo,\l)$
(\cite{A1}) a small neighborhood of the torus in $M^4$ is
biholomorphically equivalent to a neighborhood of the zero section
in $N_LM$. In \cite{A2} Arnold asks about non-integrable version
of this result.

 \begin{prop}\po
Codimension of the set of almost complex structures, the germs of
which on the PH-curve $L\subset(M,J)$ are isomorphic to these of
the normal bundle $L\subset N_LM$, in the set of all almost
complex structures is infinity.
 \end{prop}

 \begin{proof}
For existence of such an isomorphism two conditions must fulfil.
First, by Corollary \ref{999}, the Nijenhuis tensor characteristic
distribution $\Pi^2$ should be integrable and transversal to $L$
whenever non-zero. Second, by Proposition \ref{135}, the Nijenhuis
tensor $N_J$ should be constant along the leaves of $\Pi^2$. Both
conditions are of $\op{codim}=\infty$.
 \end{proof}

 \vspace{4pt}
The two mentioned conditions are necessary, but not sufficient.

 \begin{ex}
Let $M^4=L^2\times D^2$ have coordinates $(z=x+iy,w=s+it)$. Equip
it with the almost complex structure
 \begin{equation}\label{251000}
J\p_x=a_1\p_x+(1+a_2)\p_y+b_1\p_s+b_2\p_t,\quad J\p_s=\p_t.
 \end{equation}
Then $L\times\{0\}$ is a PH-curve, if $b_i=0$ on it. Moreover, one
can achieve $a_i|_L=0$.

The integrability condition $\Pi^2=T\mathcal{F}$,
$\mathcal{F}_c=\{z=c\}$, and the requirement of the tensor $N_J$
constancy along $\mathcal{F}$ write as follows ($c_i=c_i(x,y)$):
 \begin{equation*}
 \left\{
 \begin{array}{lll}
 \ddfrac{\p a_1}{\p t}=a_1\ddfrac{\p a_1}{\p s}-
\ddfrac{1+a_1^2}{1+a_2}\ddfrac{\p a_2}{\p s},&\!\!\! &
 \ddfrac{\p b_1}{\p t}=-\ddfrac{\p b_2}{\p s}+
b_1\ddfrac{\p a_1}{\p s}+
 \ddfrac{b_2-b_1a_1}{1+a_2}\ddfrac{\p a_2}{\p s}+c_1,\\
 \ddfrac{\p a_2}{\p t}=(1+a_2)\ddfrac{\p a_1}{\p s}-
a_1\ddfrac{\p a_2}{\p s}, &\!\!\! &
 \ddfrac{\p b_2}{\p t}=\ddfrac{\p b_1}{\p s}+
b_2\ddfrac{\p a_1}{\p s}-
 \ddfrac{b_1+b_2a_1}{1+a_2}\ddfrac{\p a_2}{\p s}+c_2.
 \end{array}
 \right.
 \end{equation*}

This is a Cauchy-Kovalevskaya type system, so any analytical
initial condition $(a_i,b_i)|_{t=0}=(\a_i^0(s),\b_i^0(s))$
determines uniquely the solution. PH bundle structures correspond
to $\a_i^0=\lambda_i(x,y)$, $\b_i^0=\mu_i(x,y)+\nu_i(x,y)s$. There
are however different solutions, for example:
$a_1=-b_1=-\frac{s}{1+t}$, $a_2=-b_2=-\frac{t}{1+t}$.
 \end{ex}

Thus the answer to Arnold's question is negative. A generalization
of his theory should look differently. It will concern existence
of a PH-foliation of a $T^2$-neighborhood by cylinders. In
holomorphic situation there exists a foliation by holomorphic
cylinders, given in the representation (\ref{e1}) as
$\{w=\op{const}\}$. Does it persist if we perturb the structure
$J$ to an almost complex one?

We discuss this question in \cite{K3}. Note however that in the
complex situation transport along the leaves of the foliation
is holomorphic. When does a PH-foliation exist with
pseudoholomorphic transports?

By transports here we mean the following. Let $D^2_z$ be a
foliation by transversal PH-disks as in proposition \ref{prop:1}.
Let $\mathcal{H}$ be a PH-foliation with $L$ as a leaf. A path
between two points $z_1,z_2\in L$ determines a map $D^2_{z_1}\to
D^2_{z_2}$ of shifts along $\mathcal{H}$, called the transport.
Homotopically non-trivial loops yield the monodromy (for a
PH-foliation $\mathcal{H}$ by cylinders, one cycle has a trivial
monodromy).

The requirement of PH-transports is independent of the choice of
transversal disks family. For a generic almost complex structure
the monodromy and transports are non PH-maps of the disks $D^2_z$.

 \begin{prop}\po\label{prop:5}
Let $L\subset(M^4,J)$ be a PH-curve. Existence of a PH-foliation
$\mathcal{H}$ of its neighborhood with PH-transports is a
condition of codimension infinity.
 \end{prop}

 \begin{proof}
The requirements of PH-transports means that projection along
$\mathcal{H}$ is a PH-map. Thus by Corollary \ref{999} the
Nijenhuis tensor characteristic distribution is integrable and
tangent to $\mathcal{H}$. Also the Nijenhuis tensor should be
locally projectible along $\mathcal{H}$. These are two conditions
of $\op{codim}=\infty$.
 \end{proof}

 \vspace{4pt}
Here is another generalization of Arnold's theory of holomorphic
curves neighborhoods:

 \begin{theorem}\po
A small neighborhood $\ok_L$ of a PH-curve $L=\Sigma^2_g$ is
Kobayashi hyperbolic iff $g\ge2$. For $g=0$ the punctured
neighborhood $\ok_L\setminus L$ is not hyperbolic and for $g=1$ it
is not hyperbolically imbedded into $\ok_L$.
 \end{theorem}

We refer to \cite{Kob} for the basics of hyperbolic spaces. In
almost complex context the corresponding notions were introduced
in \cite{KO} and a non-integrable version of Brody criterion was
established. Its application together with a theorem of Lang
(\S3.6 \cite{Kob}) and compactness from \cite{Gro} yield the above
statement.

\section{\hps Deformations of PH-curves}
\label{sec6}

 \hspace{13.5pt}
In this section we continue to study PH-curves
$L\simeq(\Sigma^2_g,j)$. Let $\mathcal{X}=C^\infty(\Sigma,M;A)$ be
the space of all smooth maps $u:\Sigma^2_g\to M^{2n}$ representing
a fixed homology class $A\in H_2(M)$ and
$\varrho:\mathcal{E}\to\mathcal{X}$ be the bundle with the fiber
$\varrho^{-1}(u)=\mathcal{E}_u=\Omega^{0,1}(u^*TM)$ being the
space of anti-linear maps $T\Sigma\to TM$ over $u$. For Fredholm
theory these spaces should be completed to appropriate functional
spaces (\cite{MS}), whose precise choice is not crucial due to
elliptic regularity. But we will not specify them, because it is
irrelevant for our geometric approach.

PH-curves $L=\op{Im}[u:(\Sigma^2_g,j)\to(M,J)]$ in the class $A$
are zeros of the section $\bar\p_J=\frac12(\1+J\circ j^*)\circ
d:\mathcal{X}\to\mathcal{E}$ and their union forms the moduli
space $\mathcal{M}(A,J)=\bar\p_J^{-1}(0)$. To study regularity of
a point $u\in\mathcal{M}(A,J)$ Gromov \cite{Gro} considers the
linearization
$D_u=D\bar\p_J:C^\infty(u^*TM)\to\Omega^{0,1}(u^*TM)$. This
Gromov's operator can be explicitly written (\cite{MS,IS}) as
 \begin{equation}\label{MIS}
D_u(v)=\bar\p_{u,J}(v)+\tfrac14N_J(v,\p_J(u)).
 \end{equation}
The operator descends to the normal bundle in virtue of the
following diagram:
 $$
 \begin{CD}
0\to\ @. C^\infty(T\Sigma) @>{du}>> C^\infty(u^*TM)
@>{\text{proj}}>>
C^\infty(u^*TM)/C^\infty(T\Sigma) @.\ \to 0 \\
@. @V{\bar\p}VV @V{D_u}VV @V{\check D_u}VV @. \\
0\to\ @. \Omega^{0,1}(T\Sigma) @>{du}>> \Omega^{0,1}(u^*TM)
@>{\text{proj}}>> \Omega^{0,1}(u^*TM)/\Omega^{0,1}(T\Sigma) @.\
 \to 0.
 \end{CD}
 $$

As before we consider only regular PH-curves $L=u(\Sigma)$
(singularities may enlarge the sheaf of holomorphic sections of
the normal bundle, see \cite{IS}), in which case
$C^\infty(u^*TM)/C^\infty(T\Sigma)=C^\infty(N_LM)$ and similar for
$\Omega^{0,1}$.

 \begin{prop}\po\label{D_u}
The Gromov operator $\check
D_u:C^\infty(N_LM)\to\Omega^{0,1}(N_LM)$ coincides with the
Cauchy-Riemann operator $\bar\p_{\check J}$ of the NB-II structure
$\check J$.
 \end{prop}

Of course, an indication of this result is Proposition \ref{proo}.

 \begin{proof}
This follows from Theorem \ref{th:1} because the operator
$\bar\p_{\check J}$ with the PH-bundle structure from formula
(\ref{222}) coincides with the expression (\ref{MIS}).
 \end{proof}

 \vspace{4pt}
Now we introduce the Dolbeault cohomology groups
$H^{0,0}_{\bar\p_{\check J}}(N_L^{II}M)=\op{Ker}(\check D_u)$ and
$H^{0,1}_{\bar\p_{\check J}}(N_L^{II}M)=\op{CoKer}(\check D_u)$
(of course, Sobolev spaces are needed to insure via Fredholm
property that the dimensions are finite). Vanishing of the former
is equivalent to non-existence of deformations for the PH-curve
$L$, while vanishing of the latter means transversality of
$\bar\p_J$ to the zero section of $\varrho$ at $u$, whence $u$ is
a regular point of the moduli space $\mathcal{M}(A,J)$.

Consider the case $\dim M=4$, where $N_L^IM$ is a holomorphic line
bundle. The following statement is essentially contained in
\cite{Gro,HLS,IS} (the two statements below are equivalent via
Kodaira-Serre duality).

 \begin{theorem}\po
If $c_1(N_L^IM)<0$, then $H^{0,0}_{\bar\p_{\check
J}}(N_L^{II}M)=0$. If $c_1(N_L^IM)>2g-2$, then
$H^{0,1}_{\bar\p_{\check J}}(N_L^{II}M)=0$. \qed
 \end{theorem}

For higher-dimensional $M$ in the case $g=0$ one proceeds as
follows: By Grothendieck's theorem a holomorphic line bundle over
$S^2$ splits into line bundles $N_L^IM=\oplus\mathfrak{L}_i$ and
then one gets the vanishing theorem requiring the corresponding
inequality for the Chern class of each line bundle
$\mathfrak{L}_i$.

Now for the rest of the section we study a particular interesting
case of an elliptic PH curve ($g=1$) in four-dimensional manifold
$M^4$ and its deformation. We wish to get a non-deformation
criterion, which is based on the whole structure of $N_L^{II}M$,
not only of $N_L^IM$.

Let self-intersection number of the curve be $L\cdot L=p$. If
$p<0$, the curve is not deformed by the positivity of
intersections (\cite{M1}). For $p>0$ the virtual moduli space has
positive dimension (by the index computation for the linearized
Cauchy-Riemann operator).

Consider now topologically trivial normal bundles, $p=0$, when the
elliptic curves are generically discrete and persistent under a
small perturbation of the structure $J$ (this case was studied in
\cite{Ku} and the number of non-parametrized PH-tori was
estimated). We will formulate an explicit sufficient condition of
non-deformation and persistence.

Let $\pi:(E,J)\to T^2$ be a PH line bundle. Due to corollary
\ref{coro} there exist coordinates $(z=x+iy,w=s+it)$ on $E$, with
the gluing rule (\ref{e1}), such that
 $$
\left\{
 \begin{array}{ll}
J\p_x=\hspace{8pt}\p_y+s\cdot\x-t\cdot J\x, &
J\p_s=\hspace{8pt}\p_t,
\\
J\p_y=-\p_x-s\cdot J\x-t\cdot\x, & J\p_t=-\p_s.
 \end{array}
\right.
 $$
The vertical vector field $\x=\frac12JN_J(\p_s,\p_x)$ can be
decomposed $\x=2\b_1\p_t-2\b_2\p_s$ with $\b=\b(z)$,
$\b=\b_1+i\b_2$.

Every PH-curve in $E$, homologous to the zero section $T^2$, is of
the form $L=f(T^2)$ for some section $f\in C^\infty(\pi)$. This
follows from the properties of PH-maps $\pi:L\to T^2$ of degree 1
(also from positivity of intersections for the
spheres-compactification of the fibers). Let us deduce the
equation for $f=f_1+if_2$. The tangent bundle to the curve
$w=f(z)$ is spanned by
 \begin{equation}\label{444}
\e_1=\p_x+\p_x(f_1)\p_s+\p_x(f_2)\p_t,\quad
\e_2=\p_y+\p_y(f_1)\p_s+\p_y(f_2)\p_t.
 \end{equation}
The pseudoholomorphicity condition $J\e_1=\e_2$ along $f(T^2)$ is
equivalent to the equation
 \begin{equation}
f_{\bar z}+\b\bar f=0. \label{e:f}
 \end{equation}

Below we use the normalization of \S\ref{sec5} for the pair
$(\oo,\l)$ from (\ref{e1}), characterizing the holomorphic line
bundle NB-I over an elliptic curve $L=T^2$.

 \begin{prop}\po\label{chi}
Let $J$ be a PH line bundle structure and the corresponding
complex structure $J_0$ from {\rm (\ref{222})} have the multiplier
$\l$ (normalized as in \S\ref{sec5}). Determine the function
$\Lambda\in C^\infty(T^2)$ from the equation
$\frac12JN_J(\p_w,\p_z)=\Lambda\p_{\bar w}$. Let $\l\ne1$ if
$\La\equiv0$ and if $\La\not\equiv0$ assume the inequality:
 \begin{equation}\label{rrt}
|\p_{\bar z}\La|\le(1-\ve)|\La|^2-|\tau\La|,\quad
\tau=\ln|\l|/\op{Im}\oo,
 \end{equation}
for some $\ve>0$. Then the zero section $T^2$ is the only PH-torus
in $E$.
 \label{p-pro}
 \end{prop}
Notice that if $\La\ne0$, the inequality can be achieved via a
simple rescaling. Its meaning is then that the structure $J$ is
far from being integrable ($J_0$).
 \vspace{4pt}

 \begin{proof}
We have $\Lambda=-2i\bar\b$ because
  \begin{equation}\label{beom}
 \begin{array}{ll}
\ddfrac12JN_J(\p_w,\p_z)=-2(\b_2+i\b_1)\p_{\bar w},
&\ddfrac12JN_J(\p_{\bar w},\p_z)=0,\\
\ddfrac12JN_J(\p_w,\p_{\bar z})=0, &\ddfrac12JN_J(\p_{\bar
w},\p_{\bar z})=-2(\b_2-i\b_1)\p_w,
 \end{array}
  \end{equation}
Let us show that equation (\ref{e:f}) has no nonzero solutions.

Our torus neighborhood is the product of the cylinder ${\cal
C}^2=\{z\in\C\,|\,\op{Im}z\in [0,\op{Im}\oo)\}/2\pi\Z$ and $\C(w)$
glued by the rule $(z,w)\mapsto(z+\oo,\l w)$. The boundary
$\p{\cal C}^2$ consists of the circle $S^1=\R/2\pi\Z$ and its
$\oo$-shift.

Introduce the real-valued linear function $\z=i\tau\frac12(z-\bar
z)$ and note that the function $h=e^{2\z}$ satisfies:
$h(z+2\pi)=h(z)$, $h(z+\oo)=h(z)/|\l|^2$. So using formula
(\ref{e:f}) and its consequence $f_{\bar zz}=-\b_z\bar f+|\b|^2f$
we get:
 \begin{multline*}
0= \Bigl(\frac{\l\bar\l}{\,|\l|^2}-1\Bigr)\oint_{S^1}\dfrac i2\,
e^{2\z} f_{\bar z}\bar f d\bar z=  \iint_{{\cal C}^2}\dfrac i2
d(e^{2\z} f_{\bar z}\bar f d\bar z)=\\
 \iint_{{\cal C}^2}e^{2\z}\Bigl[\,2|\b|^2|f|^2-(i\tau\b+\b_z)\bar f^2
\,\Bigr] \dfrac i2dz\we d\bar z.
 \end{multline*}
Taking the real part we deduce:
 $$
\iint_{{\cal C}^2}e^{2\z}(2|\b|^2-|\tau\b|-|\b_z|)|f|^2\,dx\we
dy\le0.
 $$
Since by the assumption $2(1-\ve)|\b|^2-|\tau\b|\ge|\b_z|$ for
some positive $\ve$, we should have $f=0$ or $\b=0$. If
$\b\equiv0$, then the Fourier decomposition of the $2\pi$-periodic
holomorphic function $f$ and the condition $f(z+\oo)=\l f(z)$ for
$e^{-\op{Im}\oo}<|\l|\le1$, $\l\ne1$, imply $f\equiv0$. If $\b$
vanishes only on a domain $D\subset T^2$, then $f$ is holomorphic
in $D$ and vanishes in $T^2\setminus D$, whence $f\equiv0$.

So there are no PH-tori $\ti T^2$, homologous to the zero section,
with $f\not\equiv0$. If the homology class of $\ti T^2$ is a
multiple of the zero section $[\ti T^2]=k[T^2]$ a $k$-finite
covering finishes the proof.
 \end{proof}

 \begin{prop}\po
The linearized equation for close PH-tori can be written as
 \begin{equation}\label{e:ff}
f_{\bar z}+\a f+\b\bar f=0.
 \end{equation}
The function $\a=0$ for the normal coordinate $w$ on $N_{T^2}M$,
equipped with NB-I complex structure and with the gluing
rule (\ref{e1}). Alternatively $\a=\op{const}$ for a global
well-defined coordinate $w$.
\end{prop}

 \begin{proof}
Since we are interested in the linearized equation, which is
determined by 1-jet of $J$, we can use the normal form given by
Theorem \ref{th:3} (we simplify it for dimension 4):
$J\stackrel{\text{\bf..}}=J_0+\frac12JN_J(r,\cdot)$. We write the
complex structure $J_0$ in coordinates $(z,w)$ of $\ok(T^2)$ with
the gluing rule (\ref{e1}) ($p=0$): $J_0\p_z=i\p_z$,
$J_0\p_w=i\p_w$. Note that in these coordinates $r=w\p_w+\bar
w\p_{\bar w}$.

The most general form of the Nijenhuis tensor along $T^2$ is the
following:
$-\frac12JN_J(\p_z,\p_w)\stackrel{\text{\bf.}}=a\p_{\bar
z}+b\p_{\bar w}$, where $a=a(z,\bar z)$, $b=b(z,\bar z)$ are
smooth functions on $T^2$. Then we obtain ($\{z=\op{const}\}$ is
assumed a PH-foliation, as in Proposition \ref{prop:1}):
 $$
J\p_z\stackrel{\text{\bf..}}=i\p_z+aw\p_{\bar z}+bw\p_{\bar w},
\quad J\p_w=i\p_w.
 $$
If $w=f(z,\bar z)$ is a surface, then $\eta=\p_z+f_z\p_w+\bar
f_z\p_{\bar w}$ and $\bar\eta=\p_{\bar z}+f_{\bar z}\p_w+\bar
f_{\bar z}\p_{\bar w}$ span a complexified tangent plane to its
graph.

Thus $w=f(z,\bar z)$ is a PH-curve iff
$J\eta-i\eta-aw\bar\eta\stackrel{\text{\bf.}}=0$ ($w=f$ is a
function of the first order of smallness on $T^2$, so we disregard
$wf_{\bar z}$), which is equivalent to the equation $f_{\bar
z}=\b\bar f$ with $\b=\frac i2\bar b$. The first statement is
proved.

We obtain the second statement, introduce a global coordinate by
the change $w\mapsto w\cdot\exp\Bigl(\dfrac{z-\bar z}{\bar\oo-\oo}
\ln\l\Bigr)$, which yields equation (\ref{e:ff}) with $\a=\dfrac
i2\dfrac{\ln\l}{\op{Im}\oo}$.
 \end{proof}

 \begin{rk}\po
Equation (\ref{e:ff}) with $\a=0$, $\b=\const$ was considered by
Moser \cite{Mo}. The proposition proves a remark on p.\,430 that
{\em "the linearized equation can be brought into form
(\ref{e:ff}) with $\a=\op{const}$"\/}.
 \end{rk}

 \begin{theorem}\po
Let the normal bundle of a PH-curve $T^2\subset(M^4,J)$ be
topologically trivial and its NB-II structure be described by the
function $\La$, as in Proposition \ref{chi}, satisfying
inequality (\ref{rrt}). Then the curve is isolated and persistent
under small perturbations of $J$.
 \end{theorem}

 \begin{proof}
This follows from Propositions \ref{proo} and \ref{chi}.
Alternatively, since index of the linearized Cauchy-Riemann
operator $P(f)=f_{\bar z}+af+b\bar f$, $f\in C^\infty(T^2,\C)$, is
zero, the required properties follow from non-existence of
non-zero solutions of the equation $P(f)=0$.
 \end{proof}

 \vspace{4pt}
Certainly a big perturbation of $J$ can destroy the properties.
Another criterion of deformations non-existence with an additional
requirement of complex transports is given by
Proposition~\ref{prop:5}.

\appendix
\section{\hps Minimal almost complex connections}
\label{a_A}

 \hspace{13.5pt}
In this appendix we prove a theorem, which is basically due to
Lichnerowicz. Our proof, however, differs from the original one
(\cite{L}).

Recall that a linear connection on an almost complex manifold
$(M,J)$ can always be taken $J$-linear. In fact, for any
connection $\na$ we can define
 $$
\hat\na_X=\frac12(\na_X-J\na_XJ).
 $$
One easily checks that $\hat\na$ is a linear connection satisfying
$\hat\na JY=J\hat\na Y$.

Also let us recall that every tensor uniquely decomposes into its
$J$-linear and anti-linear parts. For instance if $T$ is a
$(2,1)$-tensor, it has the decomposition
 $$
T=T^{++}+T^{+-}+T^{-+}+T^{--},\qquad\text{ where }
 $$
 $$
T^{\ve_1\ve_2}(JX,Y)=\ve_1JT^{\ve_1\ve_2}(X,Y),\quad
T^{\ve_1\ve_2}(X,JY)=\ve_2JT^{\ve_1\ve_2}(X,Y);
 $$
 $$
T^{\ve_1\ve_2}(X,Y)=\frac14\big[T(X,Y)-\ve_1JT(JX,Y)
-\ve_2JT(X,JY)-\ve_1\ve_2T(JX,JY)\big].
 $$

 \begin{theorem}\po\label{th_A}
For any almost complex connection $\na$ the totally antilinear
part of its torsion is $T_\na^{--}=\frac14N_J$. There are
connections, called minimal, for which $T_\na=\frac14N_J$. These
connections are sections of an affine bundle $\mathcal{M}_{(M,J)}$
associated with the vector bundle $S^2T^*M\ot_\C TM$ over $M$.
 \end{theorem}

 \begin{proof}
The first formula follows directly from the definitions. There are
also other formulae expressing the Nijenhuis tensor via a
covariant differentiation (see \cite{K1} for flat connections).

Consider now an almost complex connection $\na$. We can make a
gauge transformation $\na\mapsto\ti\na=\na+A$, $A\in
C^\infty(T^*M\ot(T^*M\ot_\C TM))$, with the $J$-linearity
condition imposed to keep $\ti\na$ almost complex. Then the
torsion is changed by the rule:
 $$
T_{\ti\na}=T_\na+\square(A),
 $$
where $\square=\op{alt}:\ot^2T^*M\ot TM\to\La^2T^*M\ot TM$ is the
alternation operator. Introducing the decomposition
 $$
A=A^++A^-,\qquad A^\ve(X,Y)=\frac12\big[A(X,Y)-\ve JA(JX,Y)\big],
 $$
we compute the components of
$\square(A)=\sum_{\ve_i=\pm}\square^{\ve_1\ve_2}(A)$:
 $$
\square^{++}(A)=A^+-A^+\tau,\ \square^{+-}(A)=-A^-\tau,\
\square^{-+}(A)=A^-,\ \square^{--}(A)=0,
 $$
where $\tau(X,Y)=(Y,X)$. We can make all the components of the
torsion vanishing, save for $T_\nabla^{--}$, using the graded
commutation relations
$T^{\ve_1\ve_2}_\na\circ\tau=-T^{\ve_2\ve_1}_\na$. Actually we get
a minimal connection $\ti\na$ with the gauge
 $$
A=-\frac12T^{++}_\na-T^{-+}_\na.
 $$
This proves the second part of the statement.

The last one follows from the above formulae for
$\square^{\ve_1\ve_2}(A)$: The gauge transformation
$\na\mapsto\na+A$ does not change the minimality iff $A^+$ is
symmetric and $A^-=0$, i.e.\ $A\in C^\infty(S^2T^*M\ot_\C TM)$.
 \end{proof}

 \begin{prop}\po
Let $\na\in\mathcal{M}_{(M,J)}$ be a minimal connection. Then
 $$
4\mathfrak{S}\{R_\na(X,Y)Z\}=\mathfrak{S}\{N_J(N_J(X,Y),Z)\}+
\mathfrak{S}\{(\na_XN_J)(Y,Z)\},
 $$
where $\mathfrak{S}$ denotes the cyclic sum.
 \end{prop}

 \begin{proof}
This is a direct corollary of the first Bianchi's identity.
 \end{proof}

 \begin{rk}\po
Thus the field of the Nijenhuis tensors $N_J\in
C^\infty(\La^2T^*M\ot_{\bar\C}TM)$ on a manifold $M$ is not
arbitrary. For a general position tensor $N_J$ this follows also
from a result of \cite{K1}: Such a tensor field $N$ restores the
structure $\pm J$, which in turn determines $N_J$ and we obtain
the constraint $N=N_J$.
 \end{rk}

The formula of the proposition involves the curvature $R_\na$, but
neither it, nor even its anti-linear part $R_\na^{--}$ is
independent of $\na\in\mathcal{M}_M$. However we have:

 \begin{prop}\po
The operator $\square_{X\we Y}Z=\na_{N_J(X,Y)}Z-4R^{--}_\na(X,Y)Z$
is independent of $\na\in\mathcal{M}_M$, tensorial in $X,Y$ and is
an $N_J$-twisted differentiation in $Z$: $\square_{X\we
Y}(fZ)=(N_J(X,Y)f)Z+f\square_{X\we Y}Z$.
 \end{prop}

 \begin{proof}
In fact, if $\ti\na=\na+A\in\mathcal{M}_M$ is another minimal
connection, then $A\in C^\infty(S^2T^*M\ot_\C TM)$ and we
calculate: $R_{\ti\na}^{--}=R_\na^{--}+\frac14A_{N_J}$.
 \end{proof}

 \vspace{4pt}
Using $\square$ we can multiply invariants of an almost complex
structure. For instance, $\square N_J\in C^\infty(\La^2_\C
T^*M\ot(\La^2T^*M\ot_{\bar\C}TM))$. There are other ways to get
invariants -- by prolongation-projection method and via the
Fr\"olicher-Nijenhuis bracket (\cite{K1}), but they are different
from this differentiation.

\section{\hps Normal bundle in other geometries.}
\label{a_B}

 \hspace{13.5pt}
The proposed construction of the tangent and normal bundle
structures is more general and can be carried out for other
geometric structures (the submanifold $L\subset M$ should allow
restriction of the structure). One chooses a Cartan connection
$\na$ (i.e.\ preserving the structure) on $M$ with a kind of
minimality: Its torsion should be equal to the corresponding
structural (Weyl) tensor (\cite{St,KL}), which is realized (we
consider, for simplicity, the case of the first order structures)
via a splitting $\z$ of the exact sequence
 $$
0\to{\frak g}^{(1)}\stackrel{i}\longrightarrow{\frak g}\ot
T_x^*M\stackrel{\d}\longrightarrow\La^2T_x^*M\ot T_xM
\stackrel{\z}{\stackrel{\dashleftarrow}\longrightarrow}
H^{0,2}({\frak g})\to0.
 $$
Here the last term is the Spencer $\d$-cohomology group (space of
structural functions), ${\frak g}\subset\op{aut}(T_xM)$ is the
symbol of the geometric structure and ${\frak g}^{(1)}=\op{Ker}\d$
is its prolongation (\cite{KL}). The freedom in a choice of $\na$
is thus reduced to ${\frak g}^{(1)}$.

For an almost complex structure $J$: ${\frak
g}=\op{gl}_\C(T_xM,J)$ and the prolongation is ${\frak g}^{(1)}=
S^2T^*_xM\ot_\C T_xM$, cf.\ Theorem \ref{th_A}.

For a symplectic structure $\Omega$: ${\frak
g}=\op{sp}(T_xM,\Omega)$, ${\frak g}^{(1)}=S^3T^*_xM$ and a
canonical normal bundle structure $\hat\Omega$ appears. By the
symplectic neighborhood theorem (\cite{W}) it is completely
determined by the restriction $\Omega_L$ and the isomorphism class
of the normal bundle with fiber-wise symplectic structure, usually
called "symplectic normal bundle" $(\nu_L,\Omega)$.

For a Riemannian structure $g$: ${\frak g}=\op{so}(T_xM,g)$,
${\frak g}^{(1)}=0$. $\nabla$ is the Levi-Civita connection. It
splits the normal bundle $N_LM$ and leads to the normal bundle
structure $\hat g$. Another approach to $\hat g$ is similar to
(\ref{complPH}): One constructs a normal foliation $\mathcal{W}$
around $L$ via geodesics $\g\subset\mathcal{W}$ in all normal
directions $T_x^\perp L$, $x\in L$, and applies the dilatations
$R_t$ along geodesics.

Defining in this way the normal structure on $N_LM$ we obtain two
structural tensors on $L$: One original on  $L\subset M$ and the
other from the normal bundle on the zero section $L\subset N_LM$.
There are relations between these tensors. For an almost complex
structure $J$ we described them in \S\ref{sec3}.

Consider now a Riemannian metric $g$. Our structural tensors are:
Riemannian curvature $R_g$ along $L$ and the normal bundle
curvature $R_{\hat g}$ at zero section.

To describe the relations consider the curvature of the normal
bundle $R^\perp$. It is the curvature tensor of the normal
connection $\nabla^\perp$, given by the orthogonal decomposition
in $TM|_L=TL\oplus N_LM$, $R=R^{\|}+R^\perp$. Note that
$R^\perp(X,Y)=R_{\hat g}(X,Y)$ for $X,Y\in TL$ and the left-hand
side is not defined for others $X,Y$.

Let ${\rm II}:TL\ot TL\to N_LM$ be the second quadratic form of
$L$ and $A:TL\ot N_LM\to TL$ be the shape (Peterson) operator
given by $g(A(X,V),Y)=g({\rm II}(X,Y),V)$, $X,Y\in TL$, $V\in
N_LM$. The Ricci equation reads:
 $$
[R_g(X,Y)V]_{\perp}=R^\perp(X,Y)V+{\rm II}(X,A(Y,V))-{\rm
II}(Y,A(X,V)),
 $$
where $X,Y,Z\in TL$, $V\in N_LM$.

In particular when $L$ is totally geodesic ${\rm II}=0$ and $A=0$,
so that the equation mean $R_g(X,Y)=R_{\hat g}(X,Y)$ for
$X,Y\in TL$ at the points of $L$.

Similar calculations occurs in other geometries, like projective
or conformal, they can be deduced from the basic structure
equations \cite{N} of these geometries.

 \vspace{5pt}
%
%
%

\ {\hbox to 10.5cm{ \hrulefill }}

{\small \textsc{Department of Mathematics and Statistics,
University of Troms\o, Troms\o, 90-37, Norway;}} \quad {\small
kruglikov\verb"@"math.uit.no}.

\end{document}